\documentclass[fleqn]{mat01}
\usepackage{times,mathtimy,amssymb}
\begin{document}

\setcounter{page}{189}
\firstpage{189}

\def\bzeta{\mbox{$\zeta$\!\!\!\!\rule[7pt]{4.5pt}{.3pt}\,}}

\newcommand{\field}[1]{\ensuremath{\mathbb{#1}}}
\newcommand{\CC}{\field{C}}
\newcommand{\HH}{\field{H}}
\newcommand{\RR}{\field{R}}
\newcommand{\ZZ}{\field{Z}}

\def\bz{\bar{z}}
\def\bv{\bar{v}}
\def\bw{\bar{w}}
\font\sss=tir at 10pt
\def\d{\hbox{\sss{d}}}

\font\xx=msam5 at 10pt
\def\qed{\mbox{\xx{\char'03\!}}}

\title{The Weierstrass--Enneper representation using hodographic
coordinates on a minimal surface}

\markboth{Rukmini Dey}{The Weierstrass--Enneper representation of a
minimal surface}

\author{RUKMINI DEY}

\address{School of Mathematics, Harish-Chandra Research Institute,
Allahabad 211 019, India\\
\noindent E-mail: rkmn@mri.ernet.in}

\volume{113}

\mon{May}

\parts{2}

\Date{MS received 19 June 2002; revised 10 October 2002}

\begin{abstract}
In this paper we obtain the general solution to the minimal surface
equation, namely its local Weierstrass--Enneper representation, using
a system of hodographic coordinates. This is done by using the method of
solving the Born--Infeld equations by Whitham. We directly compute
conformal coordinates on the minimal surface which give the
Weierstrass--Enneper representation. From this we derive the hodographic
coordinate $\rho \in D \subset {\CC}$ and $\sigma $ its complex
conjugate which enables us to write the Weierstrass--Enneper
representation in a new way.
\end{abstract}


\maketitle

\section{Introduction}

Minimal surfaces are most commonly known as surfaces which have the
minimum area amongst all other surfaces spanning a given closed curve in
${\RR}^3$. Geometrically, the definition of a minimal surface is that
the mean curvature $H \equiv 0$ at every point of the surface. If
locally one can write the minimal surface in ${\RR}^3 $ as $(x,y,
\phi(x,y))$ the minimal surface equation $H \equiv 0$ is equivalent to
\begin{equation}\label{e1}
(1+ \phi_y^2) \phi_{xx} - 2 \phi_x \phi_y \phi_{xy} + (1+\phi^2_x)
\phi_{yy}= 0.
\end{equation}
There exists a choice of conformal coordinates $(u,v) \in \Omega \subset
{\RR}^2$ so that the surface $X(u,v)=(x(u,v), y(u,v), \phi(u,v)) \in
{\RR}^3$ satisfying the minimal surface equation is given as
follows~\cite{DHKW}:
\begin{equation*}
|X_u|^2 = |X_u|^2,\quad \langle\,X_u, X_v\,\rangle = 0,\quad \Delta_{(u,v)} X = 0.
\end{equation*}

The general solution of such an equation is called the local
Weierstrass--Enneper representation. Let $D$ be a simply connected domain
in ${\CC},$ $f$ an analytic function and $g$ a meromorphic function on
$D$. Then,
\begin{equation*}
X(\tau) = \Re \int_{\tau_0}^{\tau} \Phi \,\d\zeta,
\end{equation*}
where
\begin{equation*}
\Phi = ((1-g^2)f, i (1+g^2)f, 2 fg)
\end{equation*}
is a conformal immersion of $D$ into ${\RR}^3$ which is minimal
~\cite{O}. The immersion is regular provided that wherever $g$ has a
pole of order $m$, $f$ has a zero of at least order $2m$. Moreover, $g$
is the stereographic projection of the Gauss map.

There is a simpler representation, valid away from the umbilical points of
the surface. Let $w = g(\tau)$ and $R(w) = f(w) [{\d g}/{\d
\tau}]^{-1}$. The Gaussian curvature of the surface is $K =-4
|R(w)|^{-2}(1 + |w|^2)^{-4}$. Away from the umbilical points where $K$
vanishes, ${\d g}/{\d\tau} \neq 0$ and $R(w)$ has no pole. Thus in the
neighborhood of a nonumbilic interior point, any minimal surface can be
represented in terms of $w$ as follows~\cite{N}:
\begin{align}\label{e3a}
x(\zeta) &= x_0 + \Re \int_{\zeta_0}^{\zeta} (1-w^2) R(w) \,\d w,\\[.3pc]\label{e3b}
y(\zeta) &= y_0 + \Re \int_{\zeta_0}^{\zeta} i (1 + w^2) R(w) \,\d w,\\\label{e3c}
\phi(\zeta) &= \phi_0 + \Re \int_{\zeta_0}^{\zeta} 2 w R(w) \,\d w.
\end{align}

In this semi-expository paper we show that a system of hodographic coordinates gives us
the local Weierstrass--Enneper representation of a minimal surface. Our
method provides an easy way of calculating conformal coordinates if the
formula for the graph of the minimal surface is given, locally. This is
not given in the standard text books. Hodographic coordinates are a
natural concept in fluid mechanics where velocity fields play the role
of independent variables. It was mentioned in the context of minimal
surfaces first in \cite{K} and was used in the context of Born--Infeld
equations in \cite{W}. If one replaces $y$ by $iy$ in the Born--Infeld
equations, one obtains the minimal surface equation. Thus it is natural
to expect a general solution for the minimal surface by following
Whitham's method \cite{W} for the Born--Infeld equation. Finally we
derive the hodographic coordinates $\rho, \sigma$, complex conjugates of
each other, which enables us to write the Weierstrass representation in
a new way.

\section{Hodographic coordinates and Weierstrass--Enneper representation}

In the height representation of the minimal surface, or Monge gauge, one
writes the minimal surface equation in ${\RR}^3$ as in (1).

Introducing the complex coordinates $z = x + i y$ and $\bar{z}=x-iy$, we
define $u = \phi_{\bz}$ and $v = \phi_{z}=\bar{u}$ to reduce the second-order differential equation (1) to a pair of first-order equations:
\begin{align}
&u_{z} - v_{\bz}=0,\\
&v^2 u_{\bz} - (1 + 2uv) u_{z} + u^2 v_{z} =0.
\end{align}
The hodograph transformation interchanges the dependent and independent
variables $(z,\bz)\leftrightarrow(u,v)$. To do this we use
\begin{equation}\label{invert}
\left[ \begin{matrix} z_u & z_v\\[.3pc] \bz_u & \bz_v\end{matrix}\right]
\left[\begin{matrix} u_z & u_{\bz}\\[.3pc] v_z & v_{\bz}\end{matrix}\right]=
\left[ \begin{matrix} 1&0\\[.3pc] 0&1\end{matrix}\right]
\end{equation}
and find
\begin{align}\label{ena}
&\bz_v - z_{u} = 0,\\\label{enb}
&v^2 z_v + (1 + 2uv) \bz_v + u^2 \bz_u =0.
\end{align}
Note that we have transformed the nonlinear partial differential
equations for $u$ and $v$ into linear partial differential equations for
$z$ and $\bz$. Thus, it should come as no surprise that the minimal
surface equation has a linear representation. Following \cite{W}, we
introduce the new variables $\zeta =(\sqrt{1+4uv}-1)/(2v)$,
$\bar{\zeta}= (\sqrt{1+4uv} -1)/(2u)$ to facilitate our solution. The
inverse of this transformation is $u=\zeta/(1-\zeta{\bzeta})$.

\begin{proposition}$\left.\right.$\vspace{.5pc}

\noindent In these new coordinates the eqs $(\ref{ena})$ and $(\ref{enb})$ are
greatly simplified{\rm :}
\begin{equation}\label{ee1}
\zeta^2 \bz_{\zeta} + z_{\zeta}  =0.
\end{equation}
\end{proposition}

\begin{proof}
Using the inverse transformation $u = \zeta/(1-\zeta{\bzeta}),\ v$
its complex conjugate, and the equalities $z_{\zeta} = z_{v}v_{\zeta} +
z_{u}u_{\zeta}$ and $\bz_{\zeta} = \bz_{v}v_{\zeta} + \bz_{u}u_{\zeta}$,
we obtain $\bz_{\zeta} = (\bz_{v}{\bzeta}^{2}+\bz_{u})/(1-
|\zeta|^{2})^{2}$ and $z_{\zeta} = (z_{v}{\bzeta}^{2}+z_{u})/(1-
|\zeta|^{2})^{2}$. Then eq.~(\ref{ee1}) is equivalent to
\begin{equation}
\frac{\zeta^{2}\bz_{u} + \bz_{v}(\zeta{\bzeta})^{2}+z_{u} +
z_{v}{\bzeta}^{2}}{(1-\zeta{\bzeta})^{2}} = 0.
\end{equation}

Using the expression for $u$ and $v$ in terms of $\zeta$ and
${\bzeta}$ we rewrite (9) as
\begin{equation}
\frac{{\bzeta}^{2}z_{v} + (1+(\zeta{\bzeta})^{2})\bz_{v} +
\zeta^{2}\bz_{u}}{(1-\zeta{\bzeta})^{2}} = 0.
\end{equation}

Substracting (12) from (11) and using (8) we get zero.

\noindent Thus eqs~(\ref{ena}) and (\ref{enb}) are equivalent to eq.~(\ref{ee1}).
\hfill \qed
\end{proof}

Now taking derivative of eq.~(\ref{ee1}) with respect to ${\bzeta}$,
we obtain
\begin{equation}
\zeta^{2}\bz_{\zeta{\bzeta}} + z_{\zeta{\bzeta}} = 0.
\end{equation}

Using (13) and its complex conjugate we immediately obtain
$\bz_{\zeta{\bzeta}} = z_{\zeta{\bzeta}} = 0$ from which it
follows that
\begin{equation}
\bz = \bz_{0} + F(\zeta) + G({\bzeta}).
\end{equation}
Using (10) we find
\begin{equation}
\zeta^{2} F'(\zeta)+\bar{G}'(\zeta) = 0,
\end{equation}
so
\begin{equation}\label{e5a}
\bz = \bz_{0} + F(\zeta) - \int_{{\bzeta}_0}^{{\bzeta}}
\bar{\omega}^{2}\bar{F}'(\bar{\omega})\ \d \bar{\omega}.
\end{equation}
Moreover, we have $\phi_\zeta = \phi_{\bz} \bz_\zeta + \phi_z z_\zeta =
uF'(\zeta) - v\zeta^2F'(\zeta) = \zeta F'(\zeta)$ so that
\begin{equation}\label{e5c}
\phi = \phi_0 + \int_{\zeta_0}^{\zeta} \omega F^{\prime}(\omega)\, \d
\omega + \ \int_{{\bzeta}_0}^{{\bzeta}} \bar{\omega}
\bar{F}^{\prime}(\bar{\omega})\, \d \bar{\omega}.
\end{equation}
It is straightforward to check that the coordinates $\zeta_1=\Re \zeta$
and $\zeta_2=\Im \zeta$ are isothermal so that $\vert
X_{\zeta_1}\vert^2=\vert X_{\zeta_2}\vert^2$ and $\langle\,
X_{\zeta_1},X_{\zeta_2}\,\rangle=0$.

Rewriting (\ref{e5a}) for $x$ and $y$, we have
\begin{align}\label{e3a}
x(\zeta) &= x_0 + \Re \int_{\zeta_0}^{\zeta} (1-\omega^2) F'(\omega) \,
\d \omega,\\\label{e3b}
y(\zeta) &= y_0 + \Re \int_{\zeta_0}^{\zeta} i (1
+ \omega^2) F'(\omega) \,\d \omega,\\\label{e3c}
\phi(\zeta) &= \phi_0 +
\Re \int_{\zeta_0}^{\zeta} 2\omega F'(\omega) \, \d \omega.
\end{align}

Letting $F'(\omega)=R(\omega)$, eqs~(\ref{e3a})--(\ref{e3c}) are the
Weierstrass--Enneper representation away from the umbilical points of the
surface \cite{N}. At an umbilical point, the Gaussian curvature $K$
vanishes so $\phi_{zz}\phi_{\bz\bz} - \phi_{z\bz}^2=\ u_{\bz}v_z -
u_zv_{\bz}=0$, precisely where (\ref{invert}) has no solution.

If $F^{\prime}(\zeta)\neq 0$, we can locally introduce new variables
$\rho=F(\zeta)$, $\sigma=G(\eta)$. Locally the inverse exists when the
Gaussian curvature is finite. This follows from the fact that
$K={-4}/{(|\frac{\partial \rho}{\partial \zeta}|^2 (1+|\zeta|^2)^4)}$.
If inverse exists, $(x,y, \phi)$ can be written as

\begin{align}
x &= \frac{\rho + \sigma}{2} - \frac{1}{2} \int_{\rho_0}^{\rho}
(F^{-1}(\tilde{\rho}))^2 \d \tilde{\rho} - \frac{1}{2}
\int_{\sigma_{0}}^{\sigma}(G^{-1}(\tilde{\sigma}))^2 \d \tilde{\sigma},\\[.6pc]
y &= \frac{\sigma - \rho}{2i} - \frac{1}{2i} \int_{\rho_0}^{\rho}
(F^{-1}(\tilde{\rho}))^2 \d \tilde{\rho} + \frac{1}{2i}
\int_{\sigma_0}^{\sigma}(G^{-1}(\tilde{\sigma}))^2 \d \tilde{\sigma},\\[.6pc]
\phi &= \int_{\rho_0}^{\rho} F^{-1}( \tilde{\rho}) \d \tilde{\rho} +
\int_{\sigma_0}^{\sigma} G^{-1}(\tilde{\sigma}) \d \tilde{\sigma},
\end{align}
where $\zeta = F^{-1}(\rho) = {\partial\phi}/{\partial \rho}$ and
${\bzeta} = G^{-1}(\sigma) = {\partial \phi}/{\partial \sigma}$.
Thus
\begin{align}
x &= \frac{\rho + \sigma}{2} - \frac{1}{2} \int_{\rho_0}^{\rho}
(\phi_{\tilde{\rho}})^2 \d \tilde{\rho} - \frac{1}{2}
\int_{\sigma_{0}}^{\sigma}(\phi_{\tilde{\sigma}})^2 \d \tilde{\sigma},\\[.5pc]
y &= \frac{\sigma - \rho}{2i} - \frac{1}{2i} \int_{\rho_0}^{\rho}
(\phi_{\tilde{\rho}})^2 \d \tilde{\rho} + \frac{1}{2i}
\int_{\sigma_0}^{\sigma}(\phi_{\tilde{\sigma}})^2\ \d \tilde{\sigma},\\[.5pc]
\phi &= \phi(\rho) + \phi(\sigma).
\end{align}

This decomposition is different from that of the isothermal coordinates
$\zeta$ and ${\bzeta}$. If $\rho = \rho_1 + i \rho_2$, then it can be
shown that $|X_{\rho_1}| = |X_{\rho_2}|$ and $\langle X_{\rho_1}, X_{\rho_2}\rangle =
0$. Thus $\rho_1$, $\rho_2$ are the isothermal coordinates. The $\zeta$
system and the $\rho$ system are related by a conformal map, $F(\zeta)$
and its inverse.

The geometric meaning of $\phi_{\rho}$ is as follows: The unit normal to
the surface in the $\rho$ system is given by
\begin{equation*}
N = \frac{X_{\rho_1} \times X_{\rho_2}}{| X_{\rho_1} \times X_{\rho_2}|}
= \left(\frac{2 \ \hbox{Re}\ \phi_{\rho}}{1+ |\phi_{\rho}|^2}, \frac{2 \
\hbox{Im}\ \Phi_{\rho}}{1+ |\phi_{\rho}|^2}, \frac{ (|\phi_{\rho}|^2-1)}{1+
|\phi_{\rho}|^2}\right).
\end{equation*}
Thus $\phi_{\rho}$ is the stereographic projection of the Gauss map~\cite{O}.

\section{An example}

We consider the helicoid, $\phi=\tan^{-1}(y/x)$. We have
$u={i}/{2\bz}$ and $v={-i}/{2z}$ so that

\begin{align}
\bz &= \frac{i}{2u} = \frac{i}{2}\left[\frac{1}{\zeta} -\eta\right],\\
z&=\frac{-i}{2v} = \frac{-i}{2}\left[\frac{1}{\eta}-\zeta\right],
\end{align}
where $\eta = {\bzeta}$, from which it follows that the hodographic
coordinates are $\rho=F(\zeta)={i}/{2\zeta}$ and
$\sigma=G(\eta)={-i}/{2\eta}$. The solution $\phi$ is

\begin{align}
\phi &= -\frac{i}{2} \ln \zeta + \frac{i}{2}\ \ln\ \eta\nonumber\\
&=\frac{-i}{2}\ \ln\ \left[\frac{z}{\bz}\right]\nonumber\\
&=\tan^{-1}\left(\frac{y}{x}\right).
\end{align}
Finally, note that $R(\omega) = F'(\omega) = {-i}/{2\omega^2}$ is
the standard result for the helicoid \cite{N}.

\section*{Acknowledgement}
I would like to thank Professor Randall Kamien, for suggesting
this problem to me and for useful discussions. I would like to thank Dr.
Rajesh Gopakumar and Dr. Abhijit Mukherjee for their helpful comments.

\end{document}